\documentclass[letterpaper, 10 pt, conference]{ieeeconf}  % Comment this line out if you need a4paper

\IEEEoverridecommandlockouts                              % This command is only needed if 
\usepackage{amsmath} % assumes amsmath package installed
\usepackage{tikz}
\usepackage{physics}
\usepackage{subcaption}
\usetikzlibrary{patterns}
\usetikzlibrary{decorations.pathmorphing}
\usetikzlibrary{decorations.markings}
\usetikzlibrary{arrows.meta,bending}

\usepackage[version=4]{mhchem}
\usepackage{nomencl}
\usepackage{etoolbox}
\renewcommand\nomgroup[1]{%
  \item[\bfseries
  \ifstrequal{#1}{N}{Notation}{%
  }%
  \ifstrequal{#1}{S}{Sub- \&Superscripts}{%
  }%
  \ifstrequal{#1}{G}{Geometry}{%
  }%
  \ifstrequal{#1}{M}{Math}{%
  }%
  \ifstrequal{#1}{P}{Physics}{%
  }%
]}

\makenomenclature

\title{\LARGE \bf
A first engineering principles model for dynamical simulation of cement pyro-process cyclones
}

\author{Jan Lorenz Svensen$^{1,2}$, Nicola Cantisani$^{1}$, Wilson Ricardo Leal da Silva$^{2}$,\\ Javier Pigazo Merino$^{2}$, Dinesh Sampath$^{2}$ and John Bagterp J\o rgensen$^{1}$% <-this % stops a space
\thanks{*This work was supported by Innovation Fund Denmark, Ref. 2053-00012B}% <-this % stops a space
\thanks{$^{1}$ DTU Compute, Department of Applied Mathematics and Computer Science, Technical University of Denmark, 2800 Lyngby, Denmark
        {\tt\small jlsv@dtu.dk, nicca@dtu.dk, jbjo@dtu.dk}}%
\thanks{$^{2}$ FLSmidth Cement, 2500, Valby, Denmark
        {\tt\small wld@flsmidth.com}}%
}

\begin{document}

\maketitle
\thispagestyle{empty}
\pagestyle{empty}

%%%%%%%%%%%%%%%%%%%%%%%%%%%%%%%%%%%%%%%%%%%%%%%%%%%%%%%%%%%%%%%%%%%%%%%%%%%%%%%%
\begin{abstract}
We provide a cyclone model for dynamical simulations in the pyro-process of cement production. The model is given as an index-1 differential-algebraic equation (DAE) model based on first engineering principle. Using a systematic approach, the model integrates cyclone geometry, thermo-physical aspects, stoichiometry and kinetics, mass and energy balances, and algebraic equations for volume and internal energy.
 The paper provides simulation results that fit expected dynamics. 
 The cyclone model is part of an overall model for dynamical simulations of the pyro-process in a cement plant. This model can be used in the design of control and optimization systems to improve energy efficiency and reduce \ce{CO2} emission.
\end{abstract}
%%%%%%%%%%%%%%%%%%%%%%%%%%%%%%%%%%%%%%%%%%%%%%%%%%%%%%%%%%%%%%%%%%%%%%%%%%%%%%%%
\section{Introduction}
% On the global stage,
Cement production corresponds to $8\%$ of \ce{CO2} emissions by humans globally \cite{CO2Techreport}.
The main contributor is the production of cement clinker due to calcination of \ce{CaCO3} and fuel combustion.
Technologies such as process modifications that enable the usage of alternative materials, carbon capture and utilization, process optimization, and digitalization represent the main levers to help cement plants transit towards net-zero \ce{CO2} emissions.
The development of such digitalization, control, and optimization tools requires dynamic simulation and digital twins for the cement plant.

\begin{figure}[tb]
    \centering
        \begin{tikzpicture}    
            \draw[pattern=north east lines, pattern color=blue!50] (0.2,5.5) -- (0.2,6) -- (0.8,6) -- (0.8,5.5) -- (0.6,5) -- (0.4,5)-- cycle;
            \draw[pattern=north east lines, pattern color=blue!70] (1.2,4.8) -- (1.2,5.3) -- (1.8,5.3) -- (1.8,4.8) -- (1.6,4.3) -- (1.4,4.3)-- cycle;
            \draw[pattern=north east lines, pattern color=blue!70] (0.2,4.1) -- (0.2,4.6) -- (0.8,4.6) -- (0.8,4.1) -- (0.6,3.6) -- (0.4,3.6)-- cycle;
            \draw[pattern=north east lines, pattern color=blue!70] (1.2,3.4) -- (1.2,3.9) -- (1.8,3.9) -- (1.8,3.4) -- (1.6,2.9) -- (1.4,2.9)-- cycle;
            \draw[pattern=north east lines, pattern color=blue!70] (0.2,2.7) -- (0.2,3.2) -- (0.8,3.2) -- (0.8,2.7) -- (0.6,2.2) -- (0.4,2.2)-- cycle;

            \draw  (1.4,5.5) -- (0.8,6) -- (0.8,5.8) -- (1.2,5.3)-- (1.6,5.3) -- (1.6,5.5) -- cycle;
            \draw (0.6,5) -- (0.6,4.8) -- (1.2,5.3) -- (1.2,5.1) -- (0.8,4.6)-- (0.4,4.6) -- (0.4,5);
            \draw (1.4,4.3) -- (1.4,4.1) -- (0.8,4.6) -- (0.8,4.4) -- (1.2,3.9)-- (1.6,3.9) -- (1.6,4.3);
            \draw (0.6,3.6) -- (0.6,3.4) -- (1.2,3.9) -- (1.2,3.7) -- (0.8,3.2)-- (0.4,3.2) -- (0.4,3.6);
            %\draw (1.4,2.9) -- (1.4,2.7) -- (0.8,3.2) -- (0.8,3.0) -- (1.2,2.5)-- (1.6,2.5) -- (1.6,2.9);
            \draw (1.6,2.5)  -- (0.8,3.2) -- (0.8,3.0) -- (1.2,2.5)-- (1.6,2.5);
            \draw (1.4,2.9) -- (1.4,2.7) -- (1.6,2.5) -- (1.6,2.9);
            \draw (0.5,2.8) node{\ce{Cy5}};
            \draw (1.5,3.5) node{\ce{Cy4}};
            \draw (0.5,4.2) node{\ce{Cy3}};
            \draw (1.5,4.9) node{\ce{Cy2}};
            \draw (0.5,5.6) node{\ce{Cy1}};
            
            % Calciner
            \draw[pattern=north east lines, pattern color=red!70] (1.1,1.5) -- (1.1,2.3) -- (1.2,2.5) -- (1.6,2.5) -- (1.7,2.3) -- (1.7,1.5) -- (1.6,1.3) -- (1.2,1.3)-- cycle;
            \draw (1.4,1.9) node{Ca};
            
            % kiln and cooler
            \draw[pattern=north east lines, pattern color=teal] (3.9,1+0.3) -- (3.9,1-1) -- (5,1-1) -- (5,1-0.6) -- (4.3,1-0.6) -- (4.3,1+0.3)-- cycle;
            \draw[rotate around={-5:(0,0)},pattern=north east lines, pattern color=orange] (1.5,1-0.3) -- (4,1-0.3) -- (4,1+0.2) -- (1.5,1+0.2) -- cycle;
            \draw (2.8,0.7) node{K};
            
            \draw (3.9,1+0.3) -- (3.9,1+0.2) -- (1.65,1.4) -- (1.7,1.5) -- cycle;
            \draw (1.6,1.3) -- (1.6,1.1) -- (1.55,0.575) -- (1.2,0.7) -- (0.4,2.2) -- (0.6,2.2) -- (1.2, 1) -- (1.2,1.3) -- cycle;
            \draw (4.3,0.2) node{Co};

            % fuel and gas pipes in/out
            \draw  (3.8,0.68) rectangle (4.8,0.78);
            \draw  (1.6,1.58) rectangle (2.6,1.68);
            \draw  (0.4,6) rectangle (0.6,6.3);

            \draw[->] (4.8,0.2) -- node[anchor=south west]{Clinker Output} (5.3,0.2);
            \draw[->] (5.0,2) node[anchor=south]{Fuel Input} -- (5.0,0.73)  -- (4.85,0.73);
            \draw[->] (4.8,2) -- (4.8,1.63)  -- (2.65,1.63);
            \draw[->] (0.5,6.3) -- node[anchor=south west]{Gas Output} (0.5,6.5);
            \draw[->] (2.5,5.73) node[anchor=south]{Raw Meal Input} -- (1.5,5.73)  -- (1.5,5.55);
            \draw[->] (3.0,-0.2) node[anchor=south]{Air Input} -- (4.1,-0.2)  -- (4.1,0.0);
            \draw[->] (4.1,-0.2) -- (4.3,-0.2)  -- (4.3,0.0);
            \draw[->] (4.3,-0.2) -- (4.5,-0.2)  -- (4.5,0.0);
            \draw[->] (4.5,-0.2) -- (4.7,-0.2)  -- (4.7,0.0);
            \draw[->] (4.7,-0.2) -- (4.9,-0.2)  -- (4.9,0.0);
    \end{tikzpicture}    
    \caption{
    The pyro-section for clinker production in a cement plant consists of preheating tower of cyclones ({\color{blue}Cy}), a calciner ({\color{red}Ca}), a rotary kiln ({\color{orange}K}), and a cooler ({\color{teal}Co}).}\label{fig:production}
\end{figure}
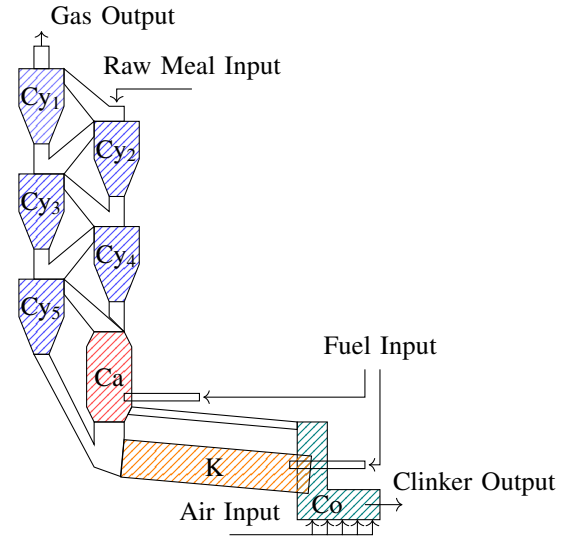

Fig. \ref{fig:production} illustrates the pyro-section of a cement plant. The pyro-section consists of a preheating tower consisting of several cyclones, a calciner, a rotary kiln, and a cooler.

In this paper, we provide a mathematical model for dynamic simulations of cyclones. 
 This model is useful for design of control and optimization systems, and is relevant for both  
 the traditional design of cement plants, as well as for modern designs, i.e., carbon capture or electrification. 

In the literature, Park et al applied CFD and finite-volume to simulate the particle and gas motion, excluding the energy and mass balances \cite{Park2020}.
Mujumdar et al proposed a steady-state model of the cyclone \cite{MUJUMDAR2007};
 the cyclone is formulated as a single-volume model (0D) with fixed efficiency and no chemical reactions.
The geometric analyses of Barth and Muschelknautz describe the cyclone efficiency and internal flows based on particle size for the steady-state case \cite{Hoffmann07}.

In contrast to the existing literature, we provide a mathematical 0D model for dynamic simulations of a cyclone, without applying CFD. 
 The model applies a systematic modeling methodology that integrates thermo-physical properties, transport phenomena, and stoichiometry and kinetics with mass and energy balances. The resulting model is a system of index-1 differential algebraic equations (DAEs). 

The mathematical simulation models of the remaining parts of the pyro-process are provided by related papers for the rotary kiln \cite{Svensen2024Kiln}, the calciner \cite{Svensen2024Calciner}, and the cooler \cite{Svensen2024Cooler}.

%Related papers provide mathematical simulation models for the rotary kiln \cite{Svensen2024Kiln}, the calciner \cite{Svensen2024Calciner}, and the cooler.

This paper is organized as follows; Section \ref{sec:Cyclone} presents the cyclone; Section \ref{sec:CycloneModel} describes the mathematical model of the cyclone; Section \ref{sec:SimulationResults} presents simulation results and conclusions are provided in Section \ref{sec:Conclusion}.

\section{The Cyclone}\label{sec:Cyclone}
In cement clinker production, the cyclone is part of the preheating tower, as shown in Fig. \ref{fig:production}. %The tower consists of two parts: 1) Risers, where a hot gas stream from below is mixed with a colder material stream from above, heating the material while flowing into a cyclone; and 2) Cyclones, where the material stream is separated from the gas stream. 
%Inside the cyclone, the materials are suspended in the gas. The outlet stream at the top of the cyclone is a gas stream with suspended materials. The outlet stream beneath the cyclone is a pure material stream, separated from the gas stream.
%While the main purpose of the cyclone is separation, the induced heating of the materials facilitates chemical reactions, e.g., calcination:
% \begin{subequations}
% \begin{alignat}{3}
% &\ce{CaCO3} + \ce{heat} \rightarrow \ce{CaO} + \ce{CO2}, \,\, &&\Delta H_r = 179.4\frac{\text{kJ}}{\text{mol}}.
% \end{alignat}
% \end{subequations}
The purpose of the tower is to facilitate efficient heat exchange between the rising hot gas stream and the falling colder material stream.
The tower consists of two parts: 1) Risers, where the falling material gets suspended in the gas stream, initiating the heat exchange, while flowing up into a cyclone; and 2) Cyclones, where the suspended material is separated from the gas stream.  The outlet stream at the top of the cyclone is a gas stream with suspended materials. The outlet stream beneath the cyclone is a stream of separated materials.
The heating of the materials facilitates chemical reactions, e.g., calcination:
\begin{subequations}
\begin{alignat}{3}
&\ce{CaCO3} + \ce{heat} \rightarrow \ce{CaO} + \ce{CO2}, \,\, &&\Delta H_r = 179.4\frac{\text{kJ}}{\text{mol}}.
\end{alignat}
\end{subequations}

\section{Dynamical Cyclone Model}\label{sec:CycloneModel}
The cyclones are described by an index-1 DAE model formulation.
The states, $x$, are the molar concentrations of each compound, $C$, and the internal energy densities of each phase, $\hat{U}$. 
We define the molar concentration vector, $C$, as mole per cyclone volume.
%The model is formulated as an index-1 DAE system, where 
The algebraic variables, $y$ are the pressure, $P$, and temperatures, $T$:
\begin{subequations}\label{eq:dyn}
\begin{align}
    \partial_tx& = f(x,y;p),\quad x=[C;\hat U],\\
    0 &= g(x,y;p),\quad y=[T;P],
\end{align}
\end{subequations}
with $p$ being a vector of system parameters.
The cyclone is described using a single-volume approach, and is formulated using a systematic modeling approach that integrates the a) geometry, b) thermo-physical properties, c) mass and energy balances, d) algebraic relations, e) transport phenomena, and f) stoichiometry and kinetics.

For the clinker compounds, we use the standard cement chemist notation, i.e.:
 \ce{(CaO)_2SiO_2} as \ce{C_2S}, \ce{(CaO)_3SiO_2} as \ce{C_3S}, \ce{(CaO)_3Al_2O_3} as \ce{C_3A} and \ce{(CaO)_4(Al_2O_3)(Fe_2O_3)} as \ce{C_4AF}, where C = \ce{CaO},   A = \ce{Al_2O_3}, S = \ce{SiO_2}, and  F = \ce{Fe_2O_3}. 

Moreover, we utilize the following assumptions:
1) the temperatures and pressure are homogeneous within the cyclone;  
2) all gasses are assumed ideal gasses;
3) the dynamics are identical across the cyclone (0D); and
4) only the five primary clinker formation reactions are included.

\subsection{Geometry}
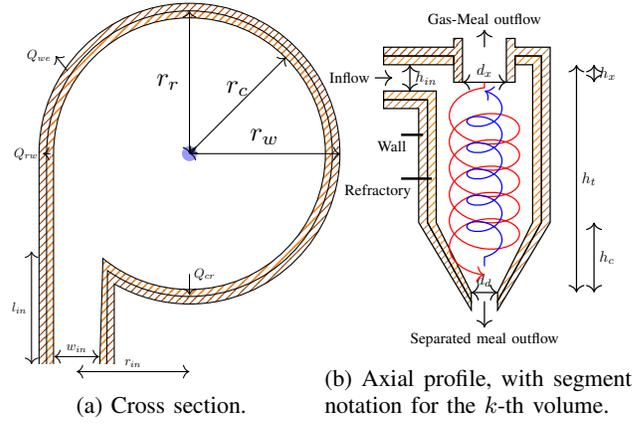
\begin{figure}
    \centering
    \begin{subfigure}[b]{0.23\textwidth}
    \centering
        \begin{tikzpicture}        
            \draw[pattern=north east lines, pattern color=brown!150] (0,-0.8) -- (0,2) arc (180:-120:2 cm) -- (1,-0.8) ; 
            \draw[pattern=north east lines, pattern color=orange] (0.1,-0.8) -- (0.1,2) arc (180:-125:1.9 cm) -- (0.9,-0.8) ;
            \draw[thick,-] (0.2,-0.8) -- (0.2,2) arc (180:-130:1.8 cm) -- (0.8,-0.8) ;
            \fill[white] (0.2,-0.8) -- (0.2,2) arc (180:-130:1.8 cm) -- (0.8,-0.8) ;
            
            \fill[blue!40!white] (2,2) circle (0.1 cm);           
    
            \draw[<->] (2,2) -- node[anchor=south]{$r_c$} (3.28,3.28) ;
            \draw[<->] (2,2) -- node[anchor=east]{$r_r$} (2.,3.9) ;
            \draw[<->] (2,2) -- node[anchor=south]{$r_w$} (4.0,2.0) ;

            \draw[->] (0.38,3.1) --  (0.22,3.3) node[anchor=east,scale=0.5]{$Q_{we}$};
            \draw[->] (0.15,2) --  (0.05,2) node[anchor=east,scale=0.5]{$Q_{rw}$};
            \draw[->] (2,0.4) -- node[anchor=south west,scale=0.5]{$Q_{cr}$} (2,0.12) ;

            \draw[<->] (0.5,-0.9)--node[anchor=south,scale=0.5]{$r_{in}$}(2,-0.9);
            \draw[<->] (0.2,-0.7)--node[anchor=south,scale=0.5]{$w_{in}$}(0.8,-0.7);
            \draw[<->] (-0.1,-0.8)--node[anchor=east,scale=0.5]{$l_{in}$}(-0.1,0.7);
        \end{tikzpicture}    
        \caption{Cross section.}
        \label{fig:cross}
    \end{subfigure} 
    \begin{subfigure}[b]{0.23\textwidth}
    \centering
    \resizebox{1\textwidth}{!}{%
    \begin{tikzpicture} 

        \draw[pattern=north east lines, pattern color=brown!150] (0.6,4.0) -- (1.4,4) -- (1.4,4.1) -- (1.5,4.1) -- (1.5,3.6) -- (1.4,3.6) -- (1.4,3.9)-- (0.6,3.9) -- cycle;
        \draw[pattern=north east lines, pattern color=orange] (0.6,3.9) -- (1.4,3.9) --  (1.4,3.8)-- (0.6,3.8) -- cycle;
        \draw[pattern=north east lines, pattern color=brown] (1,3.3) -- (0.6,3.3) -- (0.6,3.4) -- (1,3.4) --  (1.1,3.4) -- (1.1,2) -- (1.6,1.1) -- (1.6,1) -- (1,2) -- cycle;
        \draw[pattern=north east lines, pattern color=orange] (1.1,3.4) -- (0.6,3.4) -- (0.6,3.5) -- (1,3.5) --  (1.2,3.5) -- (1.2,2) -- (1.6,1.2) -- (1.6,1.1) -- (1.1,2) -- cycle;
        \draw[pattern=north east lines, pattern color=brown] (2.5,4.0) -- (2.1,4) -- (2.1,4.1) -- (2.0,4.1) -- (2.0,3.6) -- (2.1,3.6) -- (2.1,3.9) -- (2.4,3.9) -- (2.4,2) -- (1.9,1.1) -- (1.9,1) -- (2.5,2) -- cycle; 
        \draw[pattern=north east lines, pattern color=orange] (2.4,3.9) -- (2.1,3.9) -- (2.1,3.8) -- (2.3,3.8) -- (2.3,2) -- (1.9,1.2) -- (1.9,1.1) -- (2.4,2) -- cycle; 

        \draw[->] (0.5,3.65)node[anchor=east,scale=0.5]{Inflow} -- (0.7,3.65);
        \draw[->] (1.75,1.1) -- (1.75,0.8)node[anchor=north,scale=0.5]{Separated meal outflow};
        \draw[->] (1.75,4.0) -- (1.75,4.2)node[anchor=south,scale=0.5]{Gas-Meal outflow};

        \draw[<->] (2.8,3.8)--node[anchor=west,scale=0.5]{$h_t$}(2.8,1.2);
        \draw[<->] (3,2)--node[anchor=west,scale=0.5]{$h_c$}(3,1.2);
        \draw[<->] (1.6,1.2)--node[anchor=south,scale=0.5]{$d_d$}(1.9,1.2);
        \draw[<->] (1.5,3.6)--node[anchor=south,scale=0.5]{$d_x$}(2,3.6);
        \draw[<->] (0.9,3.5)--node[anchor=west,scale=0.5]{$h_{in}$}(0.9,3.8);
        \draw[<->] (3,3.6)--node[anchor=west,scale=0.5]{$h_{x}$}(3,3.8);

        \draw[thick,-] (1.05,3) -- node[anchor=north east,scale=0.5]{Wall} (0.8,3.0) ;
        \draw[thick,-] (1.15,2.5) -- node[anchor=north east,scale =0.5]{Refractory} (0.8,2.5) ;

        \draw[<-,color=red,decoration={aspect=0.61, segment length=3mm, amplitude=0.4cm,coil},decorate] (1.75,1.4) --(1.75,3.6);
        \draw[color=blue,decoration={aspect=0.61, segment length=3mm, amplitude=0.2cm,coil},decorate,arrows = {<[bend]-}] (1.75,3.5) --(1.75,1.5);        
    \end{tikzpicture}
    }
    \caption{Axial profile, with segment notation for the $k$-th volume.}
    \label{fig:axial}
    \end{subfigure}
    \caption{Geometric profiles of the cyclone.}\label{fig:profiles}
\end{figure}

Fig. \ref{fig:profiles} shows the geometry of the cyclone. It consists of a single chamber with a volume, $V_{tot}$, and a surface area, $A_c$, i.e.:
\begin{align}
    V_{tot} =\pi(r_c^2(h_t-h_c)& + \frac{h_c}{3}(r_c^2+r_x^2+r_cr_x)- r_x^2h_{x}),\\
    \begin{split}
    A_c = 2 \pi r_c(h_t-h_c)& + \pi(r_c^2-r_x^2) \\
    &  + \pi (r_c+r_d)\sqrt{(r_c-r_d)^2+h_c^2}.     
    \end{split}    
\end{align}
The collection area for the separation of the suspended solids from the gas is
\begin{small}
\begin{align}
    A_{sep} & = 2 \pi r_c(h_t-h_c) + \pi (r_c+r_2)\sqrt{(r_c-r_2)^2+h_{c,1}^2},\\
    r_2 &= r_c - \frac{h_{c,1}}{h_{c}}(r_c-r_d), \quad h_{c,1}  = h_c/2.
\end{align}
\end{small}

\subsection{Thermo-physical properties}
The thermo-physical properties of the cyclones are described for each phase by a model for the enthalpy, $H(T,P,n)$, and the volume, $V(T,P,n)$.
These models are homogeneous of order 1 in the mole vector, $n$, as
\begin{small}
 \begin{align}
    H(T,P,n) &= \sum_i n_i \left(\Delta H_{f,i}(T_0,P_0) + \int^T_{T_0}c_{p,i}(\tau)d\tau \right),
    \\
    V(T,P,n) &= \begin{cases}\sum_i n_i \left( \frac{M_i}{\rho_i} \right), \quad &\text{solid},\\ \sum_{i} n_i \left( \frac{RT}{P} \right), \quad &\text{gas}.
  \end{cases}    
\end{align}
\end{small} 
$\Delta H_{f,i}(T_0,P_0)$ is the standard enthalpy of formation at $(T_0,P_0)$. As $H$ and $V$ are homogeneous of order 1, the enthalpy and volume densities can be computed as
\begin{subequations}
\begin{align}
    \hat{H}_s &= H_s(T_m,P,C_s),\quad  &\hat{V}_s &= V_s(T_m,P,C_s),\\
    \hat{H}_g &= H_g(T_m,P,C_g),\quad &\hat{V}_g &= V_g(T_m,P,C_g),\\
    \hat{H}_r &= H_r(T_r,P,C_r),\quad &\hat{H}_w &= H_w(T_w,P,C_w).
\end{align}
\end{subequations}
The volume of each phase is given by their densities,
\begin{align}
    V_g = \hat{V}_gV_{tot},\quad  V_s = \hat{V}_sV_{tot}.
\end{align}

\subsection{Mass and Energy balance}
The cyclone's mass balances for the phase of suspended solids, $_s$, and gas, $_g$, are given for each compound $i$,
\begin{align}
    \partial_t C_{s,i} &= \frac{A_{in}N_{s,in,i}-A_{x}N_{s,x,i}-A_{sep}N_{sep,i}}{V_{tot}}  + R_i,\\
    \partial_t C_{g,i} &= \frac{A_{in}N_{g,in,i}-A_{x}N_{g,x,i}}{V_{tot}}  + R_i.
\end{align}
The energy balances are given for the mixture of solid and gas phases, $_m$, the refractory wall, $_r$, and the wall, $_w$, by
\begin{align}
    \partial_t \hat{U}_m &= \frac{1}{V_{tot}}(\Delta\tilde{H}_s + \Delta\tilde{H}_g  - Q^{cv}_{cr} - Q^{rad}_{cr}),\\
    & \Delta\tilde{H}_s = A_{in}\tilde{H}_{s,in}-A_{x}\tilde{H}_{s,x}-A_{sep}\tilde{H}_{s,sep},\\
    & \Delta\tilde{H}_g = A_{in}\tilde{H}_{g,in}-A_{x}\tilde{H}_{g,x},\\
    \partial_t \hat{U}_r &= \frac{1}{V_r}(Q^{cv}_{cr}+Q^{rad}_{cr} - Q^{cv}_{rw}),\\
    \partial_t \hat{U}_w &= \frac{1}{V_w}(Q^{cv}_{rw} - Q^{cv}_{we} - Q^{rad}_{we}).
\end{align}
$\Delta\tilde{H}_j$ is the change of enthalpy for phase $j$. The enthalpy flux is $\tilde{H}_{j,k} = H(T,P,N_{j,k})$ for the flux vector $N_{j,k}$. The heat transfer of radiation and convection is noted by $Q^{rad}$ and $Q^{cv}$.

\subsection{Algebraic equations}
The volume of the cyclone chamber is governed by the relation between the specific volume of the gas and solids,
\begin{align}
    V(T_m, P, C_g ) + V(T_m, P, C_s) = \hat{V}_{g} + \hat{V}_{s} = 1.
\end{align}
Energy conservation governs the specific energies, $\hat{U}$, relating them to temperature, pressure, and concentrations by thermo-physical properties,
\begin{small}
\begin{subequations}
\label{eq:EnergyAlgebra}
\begin{align}
    \hat{U}_m &= \hat{H}_s + \hat{H}_g - P\hat{V}_g\\
    &=H_s(T_m,P,C_s) + H_g(T_m,P,C_g) - P 
    V_g(T_m,P,C_g), \nonumber\\
    \quad \hat{U}_r &= \hat{H}_r = H_r(T_r,P,C_r), 
    \\ \hat{U}_w &= \hat{H}_w = H_w(T_w,P,C_w). 
\end{align}
\end{subequations}
\end{small}

\subsubsection{Boundary conditions}
One boundary condition of the cyclone model is the outer pressure, $P_{out}$, above the top outlet.
The second is the inflow velocity, $v_{in}$, the temperature, $T_{in}$, and the load of solid, $C_{s,in}$, and gas, $C_{g,in}$.

\subsection{Transport}
For the 0D model, the mass is transported by advection and the energy is transported by
advection, convection, and radiation.

\subsubsection{Material flux}
The mass outflow through the cyclone model consists of three fluxes describing each outflow,
\begin{align}\label{eq:fluxout}
    N_{s,x,i} &= v_{s,x}C_{s,i}, \quad N_{g,x,i} = v_{g,x}C_{g,i},\\
    N_{s,sep,i} &= v_{s,sep}C_{s,i}.
\end{align}
The inflow fluxes are given by
\begin{align}
    N_{s,in,i} &= (1-\eta_{sal})v_{in}C_{s,in,i}, \quad N_{g,in,i} = v_{in}C_{g,in,i}.
\end{align}
$\eta_{sal}$ is the separation efficiency due to saltation; particles separated immediately on entrance by flying into the wall.

\paragraph{Gas velocity}
For the gas outflow velocity, $v_{g,x}$, the cyclone pressure, $P$, is assumed representative of the entire pressure in the cyclone and located below the outlet pipe.
The velocity through the outlet pipe can be described using the turbulent Darcy-Weisbach equation with Darcy friction factor, $f_D = 0.316Re^{-\frac{1}{4}}$ \cite{Darcy-Howel},
\begin{align}\label{eq:vel}
    v_{g,x} = \Big(\frac{2}{0.316}\sqrt[4]{\frac{D_x^{5}}{\mu_m\rho_m^3}}\frac{|\Delta P|}{h_x}\Big)^{\frac{4}{7}}\text{sgn}\Big(-\frac{\Delta P}{h_x}\Big).
\end{align}
$\Delta P = P_{out}-P$ is the pressure difference. $\rho_m$ is the density of the solid-gas mixture,
\begin{align}
    \rho_m =\rho_s + \rho_g, \quad \rho_j = \sum_iM_iC_{j,i}, \ \forall j\in \{s,g\}
\end{align}
with $M$ being the molar mass. $\mu_m$ is the viscosity of the solid-gas mixture \cite{TODA2006},
\begin{align}
    \mu_m &= \mu_g\frac{1+\hat{V}_s/2}{1-2\hat{V}_s}.  \label{eq:mue}
\end{align}

\paragraph{Separation velocity}
As the particles are separated by hitting the cyclone wall, the separation flux, $N_{s,sep}$, is defined by the separation area, $A_{sep}$, and the particle radial velocity at the wall, $v_{s,sep}$. Considering the approach of Mothes and Löffler \cite{Hoffmann07}, a cylinder of equivalent volume to the cyclone is considered with radius $r_{eq}$. The separation velocity thus becomes
\begin{subequations}
\begin{align}
    v_{s,sep} &= \frac{d_m^2\Delta\rho}{18\mu_m}\frac{v_{\theta,r_{eq}}^2}{r_{eq}}, \quad r_{eq} = \sqrt{\frac{V_{tot}}{\pi h_t}},\\
    \Delta\rho &= \rho_{s,0} -\rho_g.
\end{align}
\end{subequations}
$d_m$ is the median particle diameter. $\Delta\rho$ is the difference between the solid particle density, $\rho_{s,0}$, and the gas density, $\rho_g$.
$v_{\theta,r_{eq}}$ is the tangential velocity at radius $r_{eq}$, assumed the same for gas and solids.
The tangential velocity for the radius $r$ is given by Muschelknautz as
\begin{small}
\begin{align}
    v_{\theta,r} &= \frac{ \frac{r_{c}}{r}v_{\theta,w}}{(1+\frac{f_SA_{sep}v_{\theta,w}}{2A_{in}v_{in}} \sqrt{\frac{r_c}{r}}) }, \ f_S = 0.005 (1 + 3\sqrt{c_0}).
\end{align}
\end{small}
$f_S$ is the drag friction factor. $v_{\theta,w}$ is the inlet level tangential velocity at the wall,
\begin{subequations}
\begin{alignat}{2}
 v_{\theta,w} & = \frac{r_{in}}{r_{c}\alpha}v_{in},\ \beta = \frac{w_{in}}{r_c}, \ c_0 = \frac{\sum_iM_iC_{s,in,i}}{\sum_jM_jC_{g,in,j}},\\
    \alpha & = \frac{1 - \sqrt{1 - \beta(2 - \beta)\sqrt{1 - \beta(2-\beta)\frac{1-\beta^2}{1+c_0} } }}{\beta}.
\end{alignat}
\end{subequations} 
$\alpha$ is the inlet constriction coefficient. $c_0$ is the inlet load ratio.

\paragraph{Efficiency and solid outflow velocity}
The efficiency of the cyclone, $\eta$, is defined as the ratio of separated outflow and inflow \cite{Hoffmann07},
\begin{equation}
    \eta = \eta_{sal} + \eta_{sep} = \frac{\dot m_{s,sep}}{\dot m_{s,in}}.
\end{equation}
$\eta_{sep}$ is the separation efficiency due to the internal vortex in the cyclone.
The saltation efficiency, $\eta_{sal}$, is given by
\begin{equation}
    \eta_{sal} = 1 - \min\bigg(1,\frac{c_{0L}}{c_0}\bigg).
\end{equation}
$c_{0L}$ is the cyclone loading limit,
\begin{align}\label{eq:col}
    c_{0L} &= f_c\cdot 0.025(\frac{d^*}{d_{med}})(10c_0)^k, \\
    k &= 0.15\delta + (-0.11-0.10\text{ln}(c_0))(1-\delta),\\
    \delta &= 1 \Leftrightarrow c_0 \geq 0.1.
\end{align}
$f_c$ is a correction factor. $d^*$ is the particle cut-size,
\begin{align}    
    d^* &= \sqrt{\frac{18\mu_m0.9A_{in}v_{in}}{\Delta\rho 2\pi h_i v_{\theta,r_x}^2}}.
\end{align}
$v_{\theta,r_x}$ is the tangential velocity at the outlet radius $r_x$. $h_i$ is the height of the cyclone below the outlet pipe, $ h_i  = h_t - h_{x}$.
The $0.9$ factor corresponds to an assumed $10\%$ gas flow from the inlet directly to the outlet area.

The solid mass flow at the outlet is in steady-state form defined as \cite{Hoffmann07}
\begin{align}
    \dot m_{s,x} = \dot m_{s,in} - \dot m_{s,sep} = (1 - \eta)\dot m_{s,in},
\end{align}
thus accumulation of matter is not included.
In the cyclone model, we will obtain a dynamic solid outlet flow through its flux, $N_{s,x}$, in \eqref{eq:fluxout}. We will assume the outlet velocity, $v_{s,x}$, relates to the gas outlet velocity, $v_{g,x}$,
\begin{align}\label{eq:vsout}
    v_{s,x} = f_N v_{g,x}.
\end{align}
$f_N$ is a correction factor, assumed to empirically adjust the velocity to the outlet area $A_x$.

\subsubsection{Heat convection}
The transfer of heat due to convection in the cyclone for each phase is given by
\begin{align}
    Q_{cr}^{cv} &= A_{cr}\beta_{cr}(T_m-T_r),\\
    Q_{rw}^{cv} &= A_{rw}\beta_{rw}(T_r-T_w),\\
    Q_{we}^{cv} &= A_{we}\beta_{we}(T_w-T_e),
\end{align}
where $A_{ij}$ is the in-between surface area of the mixture, refractory, wall, and environment, and $\beta_{ij}$ is the convection coefficient. 
Assuming the temperatures are located in the center of each phase, the overall convection coefficient $A\beta$ of each transfer can be formulated as
\begin{align}
    A_{cr}\beta_{cr} &= (\frac{1}{A_c\beta_m} + \frac{dx_{r_c,0.5(r_c+r_r)}}{k_r\frac{A_c+A_r}{2}})^{-1},\\
    A_{rw}\beta_{rw} &= (\frac{dx_{r_r,0.5(r_c+r_r)}}{k_rA_{r}} + \frac{dx_{0.5(r_r+r_w),r_r}}{k_w\frac{A_r+A_w}{2}})^{-1},\\
    A_{we}\beta_{we} &= \frac{k_wA_{w}}{dx_{r_w,0.5(r_r+r_w)}}.
\end{align}
$A_d$ is the surface area inside the cyclone. $dx_{i,j} = \text{ln}(\frac{r_{j}}{r_i})r_i$ is the depth for curved walls with inner and outer radius $r_i$ and $r_j$. 
The convection coefficient of the mixture, $\beta_m$, is given by \cite{Gupta2000},
\begin{align}
    \beta_m = \frac{k_m}{D_H}Nu_m.
\end{align}
$D_H = \frac{4V_{tot}}{A_d}$ is the hydraulic diameter. $k_m$ is thermal conductivity. 
$Nu_m$ is the Nusselt number,
\begin{align}
\begin{split}
    Nu_m& = 702.8 + 9.5\vdot10^{-8}\frac{v_{in}}{u_{mf}}\frac{d_c}{d_p}Re \\
    + (0.03 &+ 1.2\vdot10^{-13}\frac{v_{in}}{u_{mf}}\frac{d_c}{d_p}Re)\frac{\rho_s}{\rho_g}\frac{c_{ps}}{c_{pg}}\frac{k_s}{k_g}\frac{\Delta P_c}{0.5\rho_g v_{in}^2}.
\end{split}
\end{align}
where $d_c =2r_c$, $d_p$ is the solid particle diameter, and $u_{mf}$ is the minimum fluidization velocity (\cite{ZHOU2020Umf} suggested 0.16m/s). $Re = \rho_g v_{in} d_c/\mu_g$ is Reynolds number. $c_{p,j}$ is the specific heat capacity of phase $j$. $\Delta P_c$ is the pressure drop across the cyclone.

\subsubsection{Heat Radiation}
Radiation-driven heat transfer occurs between the solid-gas and refractory, and between wall and environment,
\begin{align}
    Q_{cr}^{rad} &= \sigma A_{p}F_{p-r}(T^4_c-T^4_r),\\
    Q_{we}^{rad} &= \sigma A_{w}(\epsilon_wT^4_w-\epsilon_eT^4_e),\\
    F_{p-r} &= \frac{1}{\frac{1}{\epsilon_p}+\frac{1}{\epsilon_r}-1}, \quad A_p=A_c\hat{V}_s.
\end{align}
$\sigma$ is Stefan-Boltzmann's constant. $\epsilon_r$ is the refractory emissivity. $\epsilon_w$ is the wall emissivity. $\epsilon_p$ is the emissivity of the particles \cite{Gupta2000}.

\paragraph{Viscosity and conductivity} 

The viscosity of a gas mixture, $\mu_g$,  and the thermal conductivity of a gas mixture, $k_g$, correlations are provided by \cite{Wilke1950} and \cite{Poling2001Book},
\begin{subequations}
\begin{align}
    \mu_g &= \sum_i\frac{x_i\mu_{g,i}}{\sum_jx_j\phi_{ij}},\quad k_g = \sum_i\frac{x_ik_{g,i}}{\sum_jx_j\phi_{ij}},\\
    \phi_{ij} &= \bigg(1+\sqrt{\frac{\mu_{g,i}}{\mu_{g,j}}}\sqrt[4]{\frac{M_j}{M_i}}\bigg)^2\bigg(2\sqrt{2}\sqrt{1+\frac{M_i}{M_j}}\bigg)^{-1}.
\end{align}
\end{subequations}
$x_i$ is the mole fraction of component $i$.
A correlation for the temperature-dependent viscosity of a pure gas is \cite{Sutherland1893}
\begin{align}
    \mu_{g,i} = \mu_0 (\frac{T}{T_0})^{\frac{3}{2}}\frac{T_0+S_\mu}{T+S_\mu}.
\end{align}
with $S_\mu$ being calibrated from two measures of viscosity.

The thermal conductivity of the solid-gas mixture is given by the serial thermal conductivity  \cite{Perry}; assuming that the solid-gas mixture can be considered as layers,
\begin{align}
    \frac{1}{k_m} = \frac{V_{g}}{V_{tot}}\frac{1}{k_g}+ \sum_i\frac{V_{s,i}}{V_{tot}}\frac{1}{k_{s,i}}.
\end{align}
The volumetric ratios describe the layer thickness.

\subsection{Stoichiometry and kinetics}
The production rates, $R$, are provided by the reaction rate vector, $r = r(T,P,C)$, and the stoichiometric matrix, $\nu$,
\begin{align}
    R = \nu^Tr.
\end{align}
The cyclone model contains the following reactions for the solid-gas mixture,
\begin{subequations}
\begin{alignat}{5}
    \text{$\#1$: }& & \ce{CaCO3} &\rightarrow \ce{CO2} + \ce{CaO}, \ & r_1,\\
    \text{$\#2$: }& & 2\ce{CaO} + \ce{SiO_2} &\rightarrow \ce{C_2S}, \quad & r_2,\\
    \text{$\#3$: }& &\ce{CaO} + \ce{C_2S}&\rightarrow \ce{C_3S}, \quad & r_3,\\
    \text{$\#4$: }& &3\ce{CaO} + \ce{Al_2O_3}&\rightarrow \ce{C_3A}, \quad & r_4,\\
    \text{$\#5$: }& &4\ce{CaO} + \ce{Al_2O_3} + \ce{Fe_2O_3}&\rightarrow \ce{C_4AF}, \quad & r_5, 
    \end{alignat}
\end{subequations}
In this paper, the rate functions, $r_j(T,P,C)$, are given by the expression
\begin{align}
    r = k(T)\prod_lC_l^{\alpha_l},\quad k_j(T) = k_{0}e^{-\frac{E_{A}}{RT}}.
\end{align}
 $k(T)$ is the Arrhenius expression. $C_l$ is the concentration (mol/L). $\alpha_l$ is the stoichiometric coefficient.
 
\begin{table}
\centering
    \caption{Reaction rate coefficients.}
    \begin{tabular}{c c|c c | c c c }
     & & $k_r$ & $E_{A}$ & $\alpha_1$ & $\alpha_2$ & $\alpha_3$\\ \hline
    Reactions&Units & $\frac{\text{kg}}{\text{m}^3s}\cdot [C]^{-\Sigma\alpha}$ & $\frac{\text{kJ}}{\text{mol}}$& 1& 1& 1\\ & && && &\\[-1.0em]\hline
    $r_1$  & $\frac{\text{kg}}{\text{m}^3s}$& $10^{8}$& $175.7$& 1& &\\[-1.0em]& && && &\\
    $r_2$ & $\frac{\text{kg}}{\text{m}^3s}$ &$10^{7}$ & $240$& 2& 1&\\[-1.0em]& && && &\\
    $r_3$ & $\frac{\text{kg}}{\text{m}^3s}$ & $10^{9}$& $420$& 1& 1&\\[-1.0em]& && && &\\
    $r_4$ & $\frac{\text{kg}}{\text{m}^3s}$ &$10^{8}$ & $310$& 3& 1&\\[-1.0em]& && && &\\
    $r_5$ & $\frac{\text{kg}}{\text{m}^3s}$ & $10^{8}$& $330$& 4& 1&1\\[-1.0em]& && && &\\\hline
    \end{tabular} 
    
    \footnotesize{The reported units and coefficients are from \cite{Mastorakos1999CFDPF}.}
       \label{tab:reaction}
\end{table}

\section{Simulations}\label{sec:SimulationResults}
\begin{figure}
    \centering
    \includegraphics[width=0.5\textwidth,trim={0.4cm, 1.1cm, 1.0cm, 1.0cm},clip]{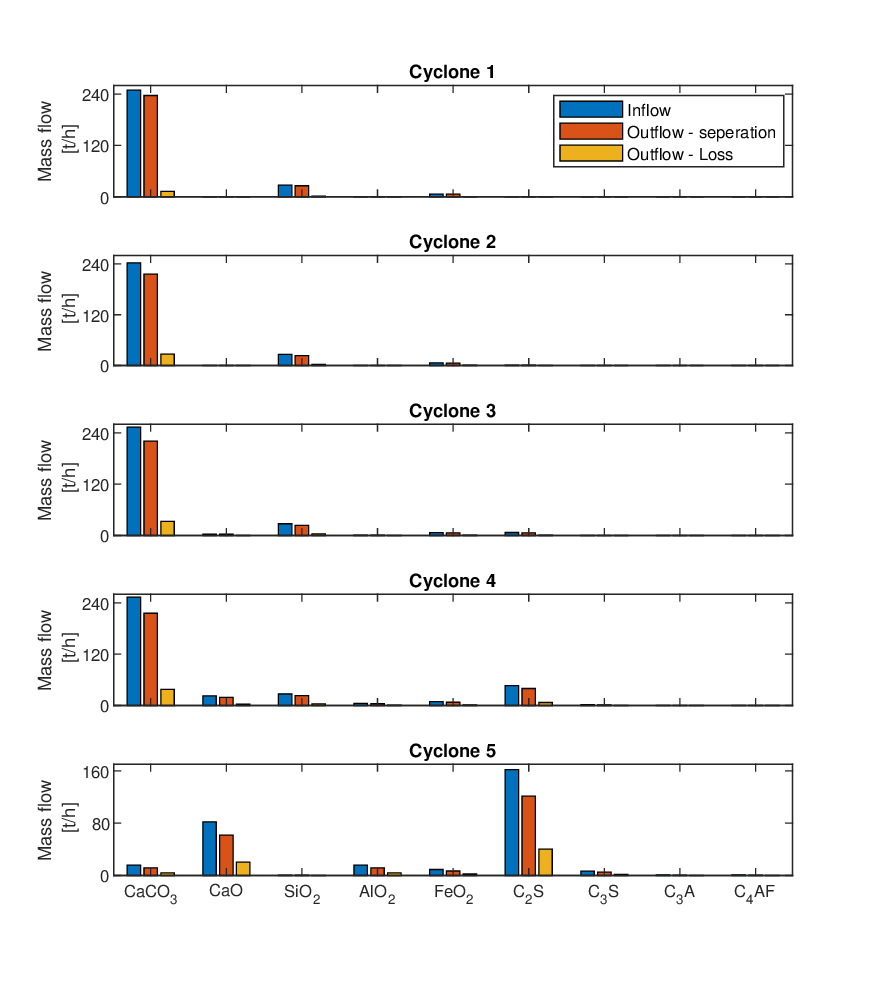}
    \caption{The solid mass flow at inlet, outlet and separation for each cyclone.}
    \label{fig:mass}
\end{figure}

To demonstrate the cyclone model, we simulate 5 scenarios each using a different cyclone for 50 hours.
Table \ref{tab:cycloneGeo} shows the dimensions of each cyclone. 
The reference scenarios are %taken from the steady-state simulator Onecalc. 
steady-state simulations from an FLSmidth Cement database.
The volumetric inflow of the solid-gas mixture, $\dot V_{in}$, in each scenario is 173.1, 223.83, 259.53, 289.59, and 309.34 $\text{m}^3/\text{s}$. Fig. \ref{fig:mass} shows the inflows of each material compound to the cyclones, and the resulting outflows.
The ambient temperature is 25$^\circ$C, and the false air inflows leaking into the cyclones are 0.95, 0.91, 0.45, 0.44, and 0.44 $\text{m}^3/\text{s}$ respectively.
Tables \ref{tab:performP} and \ref{tab:performT} give the external pressures and temperatures for each scenario.
The scenarios have the efficiencies, $\eta$, of 94.9\%, 89.0\%, 87.0\%, 85.0\%, and 75.0\%.
In each simulation, the initial concentrations of solids are uniform across compounds at 1 mol/m$^3$.
The concentration of solids in each scenario is initiated at uniform.
\begin{table}[ht]
    \centering
    \caption{ Dimensions for cyclone 1 - 5 in meters. The 
    wall thickness is 0.008 m. $w_{in}$ is $0.4r_c$.}
    \begin{tabular}{c|c c c| c c c c| c c}
         \#  & $h_t$ & $h_c$ & $h_x$ & $r_c$ & $r_r$ & $r_x$ & $r_d$  &$r_{in}$& $A_{in}$ \\ \hline
        1 &  18.3& 7.4& 3.5& 3.5 & 3.6& 1.9& 0.3& 2.8& 11.0\\
        2 &  11.4& 7.3& 3.4& 3.4 & 3.6 &2.4 &0.5 & 2.7 & 13.3\\
        3 &  11.2& 7.8& 3.4& 3.4& 3.6& 2.5& 0.5& 2.7&13.7\\
        4 &  12.0& 8.1& 3.5& 3.5& 3.7& 2.6& 0.5& 2.8&14.8\\
        5 &  12.0& 8.1& 3.5& 3.5& 3.7& 2.6& 0.5& 2.8&14.8\\ \hline
    \end{tabular}    
    \label{tab:cycloneGeo}
\end{table}

% solid inflow
   % 69.2330    0.0120    7.5110    0.0820    1.7440    0.0250    0.0010    0.0001    0.0001
   % 67.3980    0.1200    7.3090    0.1020    1.7230    0.2490    0.0110    0.0013    0.0014
   % 70.4120    0.9370    7.6140    0.2700    1.8930    1.9520    0.0820    0.0100    0.0110
   % 70.3920    6.2200    7.4710    1.3370    2.4780   12.9590    0.5470    0.0670    0.0700
   %  4.3720   22.8020    0.1380    4.3640    2.5460   44.9320    1.8970    0.2320    0.2440

% \begin{figure}
%     \centering
%     \includegraphics[width=0.5\textwidth,trim={3.3cm, 9.5cm, 4.0cm, 9.8cm},clip]{Figures/massflow.pdf}
%     \caption{The solid mass flow at inlet, outlet and separation for each cyclone.}
%     \label{fig:mass}
% \end{figure}

\begin{table}[h]
    \centering
     \caption{Steady-state pressure results. Column 4 is the mean pressure between the external pressures}
    \begin{tabular}{c|ccc|c}
    Module & $P$ & $P_{in}$ & $P_{out}$ & mean\\ \hline
    Units & bar & bar & bar & bar\\ \hline \\[-1.0em]
    Cyclone 1  &  0.9485 &   0.9529 &   0.9452  &  0.9490\\
    Cyclone 2  & 0.9584  &  0.9616  &  0.9550   & 0.9583\\
    Cyclone 3  & 0.9671  &  0.9710  &  0.9631  &  0.9671\\
    Cyclone 4  & 0.9769  &  0.9810  &  0.9729  &  0.9770\\
    Cyclone 5  & 0.9867  &  0.9906  &  0.9830  &  0.9868\\ \hline
    \end{tabular} 
    \label{tab:performP}
\end{table}

\begin{table}[ht]
    \centering
    \caption{Steady-state temperature results. column 4 is the difference between the simulated temperature and the referred outgoing temperature.}
    \begin{tabular}{c|cccc}
    Module & Temp. &ref. $T_{in}$ & ref. $T_{out}$ & Temp. diff.\\ \hline
    Units &$^\circ$C & $^\circ$C &$^\circ$C &$^\circ$C\\ \hline \\[-1.0em]
    % Cyclone 1  & 318.3454  & 321.6500&   318.9000  &  -0.5546\\
    % Cyclone 2  & 522.8109 &  526.1300&   522.6500  &   0.1609\\
    % Cyclone 3  & 674.1667 &  676.4500&   673.9300  &   0.2367\\
    % Cyclone 4  & 810.0025 &  812.7900&   809.9600 &    0.0425\\
    % Cyclone 5  & 900.5336&   903.8000 &  900.0000 &    0.5336\\ \hline
  Cyclone 1  &318.74 & 321.65 & 318.90 & -0.16\\ 
  Cyclone 2  &522.32 & 526.13 & 522.65 & -0.33\\
  Cyclone 3  &673.98 & 676.45 & 673.93 &  0.05\\
  Cyclone 4  &809.89 & 812.79 & 809.96 & -0.07 \\
  Cyclone 5  &900.65 & 903.80 & 900.00 &  0.65 \\ \hline
    \end{tabular}    
    \label{tab:performT}
\end{table}

\subsection{Calibration of steady-state}
To fit the cyclone model to the steady-state references, the models are calibrated for the following parameters.
The correction factors $f_N$ in \eqref{eq:vsout} and $f_c$ in \eqref{eq:col} are calibrated together to obtain the reported efficiency and solid density.
%$f_N$ is (6.5, 4.2, 4.85, 5.2, 6.72) and $f_c^-1$ is (22, 10.1, 8.5, 7.3, 4.2) for each cyclone.
% Based on the reported efficiency, suitable steady-state values for the scaling function $f_v$ in \eqref{eq:vsout} were found (2.642, 1.254, 1.365, 1.966, and 3.099).
% Using least-square estimation as a basis, $f_v$ is formulated as
% 6.5, 4.2, 4.85, 5.2, 6.72  
% 1/22, 1/10.1, 1/8.5, 7.3, 1/4.2
% 0.6044, 0.5579, 0.5072, , 0.3681 
%STRAUSS1987
% \begin{align}
% \begin{split}
%     f_v = L_p^T\cdot [\dot V_{in},c_0&,log(c_0), \sqrt{c_0},w_{in},r_{in},h_{in},A_{in},\\ 
%                               &r_c,r_x,r_d, h_t,h_c,h_x, p_{fall}]^T
% \end{split}
% \end{align}
%$L_p$ is found to be  [0.0340, -0.8830, -1.0372, -0.4725, 0.1158, 0.2317, -0.6167, 0.2483, 0.2896, 0.0070, 0.0206, 0.3906, -1.3856, 0.2896, -0.4260].
Table \ref{tab:perform} shows the achieved efficiency and solid density matching the reference.
For a suitable pressure drop, the steady-state pressures were calibrated to the mean of external pressures, by scaling the Darcy friction factor, $f_D$, used in \eqref{eq:vel} by 410, 815, 890, 870, and 840; corresponding to a more abrasive pipe. Table \ref{tab:performP} shows the resulting pressure.
Based on Cyclone 5, the $r_1$ reaction was tuned by a factor of 0.001 to fit the data.

\begin{table}[ht]
    \centering
    \caption{Steady-state efficiency. The saltation part of the tuned efficiencies are 0.60, 0.56, 0.51, 0.54, and 0.37 }
    \begin{tabular}{c|cc|cc|cc}
    Module & \shortstack{Efficiency\\(sim.)} &  \shortstack{Efficiency\\(ref.)} & \shortstack{$\rho_s$\\ (sim.)}& \shortstack{$\rho_s$\\ (ref.)} &\shortstack{$f_N^{-1}$\\ \quad} & \shortstack{$f_c$\\ \quad} \\ \hline
       Cyclone 1  & 94.96\% & 94.94\% & 0.499& 0.504& 22& 6.5\\
       Cyclone 2  & 89.01\% & 89.00\% & 0.378& 0.380& 10.1& 4.2\\
       Cyclone 3  & 86.94\% & 87.00\% & 0.354& 0.354& 8.5& 4.85\\
       Cyclone 4  & 85.06\% & 85.00\% & 0.380& 0.388& 7.3& 5.2 \\
       Cyclone 5  & 75.00\% & 75.00\% & 0.277& 0.277& 4.2& 6.72 \\ \hline
    \end{tabular}    
    %0.6044, 0.5579, 0.5072, 0.5399, and 0.3681 }
    %mean error in the estimated $f_v$ is $0.0166\%$.}
    \label{tab:perform}
\end{table}

% \begin{figure}
%     \centering
%     \includegraphics[width=0.5\textwidth,trim={3.5cm, 10.8cm, 4.0cm, 11.0cm},clip]{Figures/pressure_dyna.pdf}
%     \caption{The pressure profile in the dynamic simulation. $Cy_i$ is cyclone $i$.}
%     \label{fig:press}
% \end{figure}
\begin{figure}
    \centering
    \includegraphics[width=0.5\textwidth,trim={0.5cm, 0cm, 1.0cm, 0.0cm},clip]{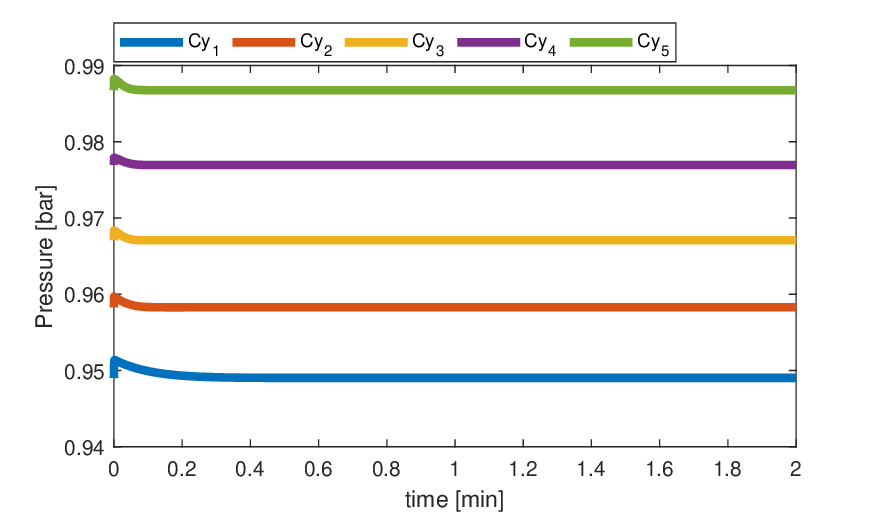}
    \caption{The pressure profile in the dynamic simulation. $Cy_i$ is cyclone $i$.}
    \label{fig:press}
\end{figure}

% \begin{figure}
%     \centering
%     \includegraphics[width=0.5\textwidth,trim={3.4cm, 9.4cm, 4.0cm, 9.5cm},clip]{Figures/temperature_dyna.pdf}
%     \caption{The temperature profile for material-gas and refractory in the dynamic simulation. $Cy_i$ is cyclone $i$.}
%     \label{fig:temp}
% \end{figure}
\begin{figure}
    \centering
    \includegraphics[width=0.5\textwidth,trim={0.3cm, 0.0cm, 1.0cm, 0cm},clip]{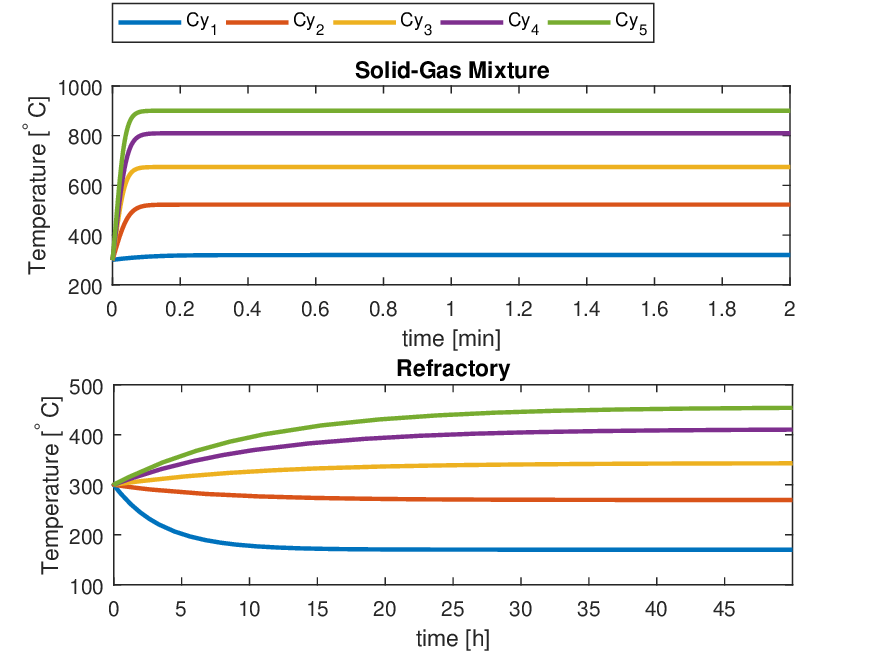}
    \caption{The temperature profile for material-gas and refractory in the dynamic simulation. $Cy_i$ is cyclone $i$.}
    \label{fig:temp}
\end{figure}

% \begin{figure}
%     \centering
%     \includegraphics[width=0.5\textwidth,trim={3.2cm, 9.0cm, 4.0cm, 8.8cm},clip]{Figures/SolidCon_dyna.pdf}
%     \caption{The concentration profiles for each material in the dynamic simulation. $Cy_i$ is cyclone $i$.}
%     \label{fig:conc}
% \end{figure}

% \begin{figure}
%     \centering
%     \includegraphics[width=0.5\textwidth,trim={3.2cm, 6.5cm, 4.0cm, 6.3cm},clip]{Figures/SolidCon_dyna1.pdf}
%     \caption{The concentration profiles for each material in the dynamic simulation. $Cy_i$ is cyclone $i$.}
%     \label{fig:conc}
% \end{figure}

% \begin{figure}
%     \centering
%     \includegraphics[width=0.5\textwidth,trim={3.2cm, 8.5cm, 4.0cm, 8.2cm},clip]{Figures/SolidCon_dyna2.pdf}
%     \caption{The concentration profiles for the main materials in the dynamic simulation. $Cy_i$ is cyclone $i$.}
%     \label{fig:conc}
% \end{figure}
\begin{figure}
    \centering
    \includegraphics[width=0.5\textwidth,trim={0.2cm, 0.4cm, 1cm, 0.0cm},clip]{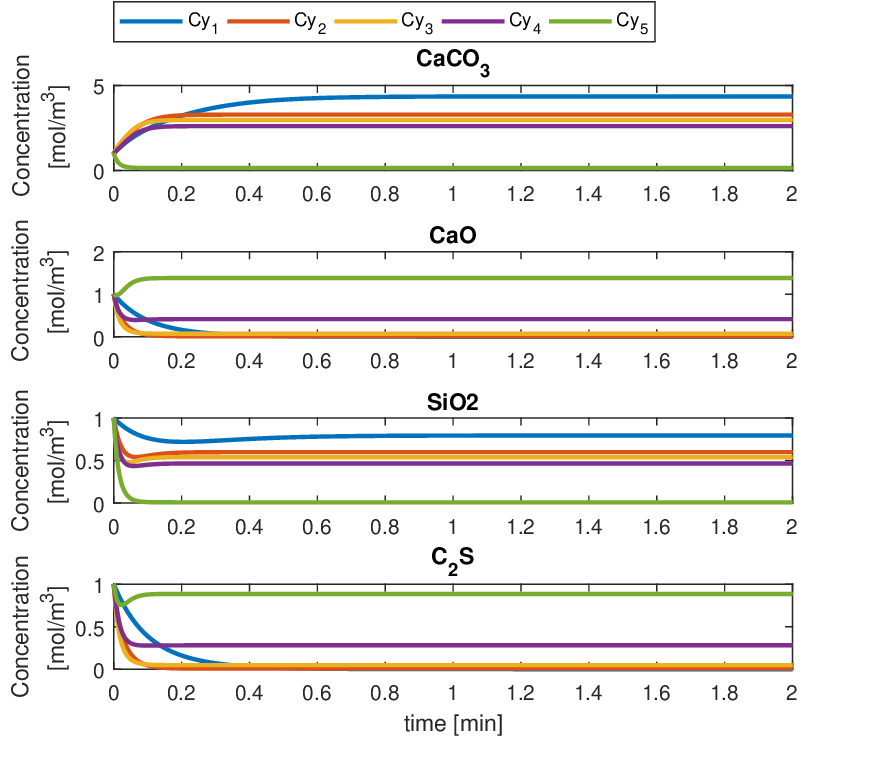}
    \caption{The concentration profiles for the main materials in the dynamic simulation. $Cy_i$ is cyclone $i$.}
    \label{fig:conc}
\end{figure}

\subsection{Simulation}
Fig. \ref{fig:press}, Fig. \ref{fig:temp}, and Fig. \ref{fig:conc} show the dynamical evolution of the cyclones. Fig. \ref{fig:press} shows the dynamics of the pressures, Fig. \ref{fig:temp} shows the temperatures of the solid-gas mixture and refractory, and Fig. \ref{fig:conc} shows the concentration of the material compounds.
The pressure, solid concentration, and solid-gas temperature all settle within 10-30 seconds, fitting to the 10 seconds on average that solid particles spend in the cyclone \cite{STRAUSS1987}. 
Cyclone 1 settles around 30 seconds, while cyclone 2-5 settles about 10-15 seconds.
Before settling down, the pressures all overshoot as the concentrations and temperatures stabilize.
For the refractory temperature, the settling time is about 10 to 30 hours.
% given the fast flow dynamics not loosing a lot of heat, and due to the absence of the sliding flow of the separation.
The steady-state solid-gas temperature in Table \ref{tab:performT} shows the cyclone temperature settles on a temperature approximately on the reported outlet temperature.

\section{Conclusion}\label{sec:Conclusion}
This paper presented a dynamic cyclone model for the case of cement clinker production in the pyro-section of cement plants. The model is an index-1 DAE model for dynamical simulations, based on a systematic approach integrating mass and energy balances with thermo-physical properties, transport aspects, reaction kinetics, and algebraic formulations for the volume and internal energy.
By calibration of a few properties, the provided simulations of the cyclone model can qualitatively match practical operations. 

For practical purposes, the proposed model can be used as part of a complete dynamic model of the pyro-process towards simplifying the application of advanced control methods.

%For dynamic simulation of the pyro-section of a cement plant, the cyclone model will be combined with models for the rotary kiln, the calciner and the cooler.
%A complete model is important for model-based design and development of control and optimization systems for the cement plant's pyro-section.

%The calciner model will be connected with models for the pre-heating cyclones, the rotary kiln, and the cooler such that the pyro-section of a cement plant can be dynamically simulated. Such a model is important for model-based design and development of control and optimization systems for the pyro-section of cement plants.

%For existing cement plants, the proposed simulation model can be used as a simplifying step in the development of advance process control systems. 
%%%%%%%%%%%%%%%%%%%%%%%%%%%%%%%%%%%%%%%%%%%%%%%%%%%%%%%%%%%%%%%%%%%%%%%%%%%%%%%%

%\section*{ACKNOWLEDGMENT}

\bibliographystyle{IEEEtran}
\bibliography{IEEEabrv,biblio}
% \newpage
\section*{APPENDIX}

\section{Properties}
Table \ref{tab:Data-Coeff-solid}-\ref{tab:heatCap} provide the parameters and physical properties used in the paper. 
Table \ref{tab:Data-Coeff-solid} and Table \ref{tab:Data-Coeff-gas} shows literature data  for the solid and gas material properties. 

The molar heat capacity of \ce{CaCO3} is described by \cite{Jacobs1981}
\begin{equation}
\begin{split}
c_p &= -184.79 + 0.32 \cdot 10^{-3}T -0.13\cdot 10^{-5}T^2 \\ & -3.69\cdot 10^{6}T^{-2} + 3883.5T^{-\frac{1}{2}} \qquad  [\frac{\text{J}}{\text{mol}\cdot \text{K}}].
\end{split}
\end{equation}
for the temperature range of 298-750 K. 

The specific heat capacities of the remaining components are computed by \cite{Svensen2024Kiln}
\begin{equation}
c_p = C_0 + C_1 T + C_2 T^2.
\end{equation}
Table \ref{tab:heatCap} reports the coefficients ($C_0$, $C_1$,$C_2$).

\begin{table}[h]
    \centering
    \caption{Material properties of the solid phase}% - perry might be profiles}
    \begin{tabular}{c|ccc}
           &\shortstack{Thermal\\ Conductivity} & Density & \shortstack{Molar \\mass}\\ \hline \\[-1.0em]
           Units    & $\frac{\text{W}}{\text{K m}}$ & $\frac{\text{g}}{\text{cm}^3}$ & $\frac{\text{g}}{\text{mol}}$ \\[-1.0em]\\ \hline 
         \ce{CaCO_3} &  2.248$^a$& 2.71$^b$  &100.09$^b$\\ 
         \ce{CaO}    & 30.1$^c$ &  3.34$^b$ &56.08$^b$\\ 
         \ce{SiO_2}   &  1.4$^{a,c}$& 2.65$^b$  &60.09$^b$\\ 
         \ce{Al_2O_3} &  12-38.5$^c$ 36$^a$& 3.99$^b$  &101.96$^b$\\ 
         \ce{Fe_2O_3} &  0.3-0.37$^c$& 5.25$^b$  &159.69$^b$\\ \hline
         \ce{C2S}     &  3.45$\pm$0.2$^d$& 3.31$^d$  &$172.24^g$\\ 
         \ce{C3S}     & 3.35$\pm$0.3$^d$ & 3.13$^d$ & 228.32$^b$\\ 
         \ce{C3A}    &  3.74$\pm$0.2$^e$& 3.04$^b$ & 270.19$^b$\\ 
         \ce{C4AF}    &  3.17$\pm$0.2$^e$& 3.7-3.9$^f$  &$485.97^g$ \\ \hline        
    \end{tabular}   
    
    \footnotesize{$^a$ from \cite{Perry}, $^b$ from \cite{CRC2022}, $^c$ from \cite{Ichim2018}, $^d$ from \cite{PhysRevApplied}, $^e$ from \cite{Du2021},\\ $^f$ from \cite{Portland}, $^g$ Computed from the above results} 
    \label{tab:Data-Coeff-solid}
\end{table}
\begin{table}[t]
    \centering
    \caption{Material properties of the gas phase}
    \def\arraystretch{1.5}
    \begin{tabular}{c|cccc}
           &\shortstack{Thermal\\ Conductivity$^a$} & \shortstack{Molar\\ mass$^a$} & Viscosity$^a$ & \shortstack{diffusion\\ Volume$^b$}\\ \hline
        Units   & $\frac{10^{-3}\text{W}}{\text{K m}}$ & $\frac{\text{g}}{\text{mol}}$ & $\mu$Pa s & cm$^3$\\[-1.5em]\\ \hline 
         \ce{CO_2}    &\shortstack{\strut 16.77 (T=300K)\\ 70.78 (T=1000K) }  & 44.01  & \shortstack{\strut 15.0 (T=300K)\\ 41.18 (T=1000K) }& 16.3\\\hline
         \ce{N_2} & \shortstack{\strut 25.97(T=300K)\\  65.36(T=1000K) }  &28.014 &  \shortstack{ \strut 17.89(T=300K)\\  41.54(T=1000K) } & 18.5\\\hline
         
         \ce{O_2}  & \shortstack{\strut 26.49(T=300K)\\  71.55(T=1000K) } &  31.998 &  \shortstack{\strut 20.65 (T=300K)\\ 49.12 (T=1000K) }&  16.3\\\hline
         
         \ce{Ar}  & \shortstack{\strut 17.84 (T=300K)\\ 43.58 (T=1000K) }& 39.948&  \shortstack{\strut  22.74(T=300K)\\  55.69(T=1000K) }&  16.2\\\hline
         
         \ce{CO}  & \shortstack{\strut 25(T=300K)\\  43.2(T=600K) } & 28.010& \shortstack{\strut 17.8(T=300K)\\  29.1(T=1000K) } &  18\\\hline
         
         \ce{C_{sus}} & - & 12.011 & - &  15.9\\\hline
         
         \ce{H_2O} & \shortstack{\strut 609.50(T=300K)\\  95.877(T=1000K) }&18.015 &  \shortstack{\strut 853.74(T=300K)\\  37.615(T=1000K) }&   13.1\\\hline
         
         \ce{H_2} & \shortstack{\strut 193.1 (T=300K)\\ 459.7 (T=1000K) }&2.016 &  \shortstack{ \strut 8.938(T=300K)\\ 20.73 (T=1000K) }& 6.12\\\hline
    \end{tabular}

    $^a$ from \cite{CRC2022}, $^b$ from \cite{Poling2001Book}
    \label{tab:Data-Coeff-gas}
\end{table}

\begin{table}[t]
    \centering
    \caption{Molar heat capacity}
    \begin{tabular}{c|c c  c| c}
            & $C_0$ & $C_1$ & $C_2$ & Temperature range\\ \hline
         Units & $\frac{\text{J}}{\text{mol}\cdot \text{K}}$& $\frac{10^{-3}\text{J}}{\text{mol}\cdot \text{K}^2}$&$ \frac{10^{-5}\text{J}}{\text{mol}\cdot \text{K}^3}$& K\\ & & & &\\[-1em]\hline
%         \ce{CaCO_3}$^a$ &  23.12 &  263.4 & -19.86& 300 - 600\\ 
         \ce{CaO}$^b$     &  71.69& -3.08  & 0.22  & 200 - 1800\\ 
         \ce{SiO_2}$^b$   & 58.91 &  5.02 & 0& 844 - 1800\\ 
         \ce{Al_2O_3}$^b$ &  233.004&-19.5913   &0.94441  & 200 - 1800\\ 
         \ce{Fe_2O_3}$^a$ & 103.9  & 0 & 0 & -\\ \hline
         \ce{C2S}$^b$     &  199.6& 0  &0  & 1650 - 1800\\ 
         \ce{C3S}$^b$     &  333.92&  -2.33&  0& 200 - 1800\\ 
         \ce{C3A}$^b$    & 260.58  & 9.58/2 & 0  &298 - 1800\\ 
         \ce{C4AF}$^b$  &  374.43& 36.4 & 0 &  298 - 1863\\ \hline
         \ce{CO_2}$^a$    & 25.98 &43.61 &-1.494 & 298 - 1500\\
         \ce{N_2}$^a$ & 27.31&5.19 &-1.553e-04 & 298 - 1500\\ 
         \ce{O_2}$^a$ & 25.82&12.63 &-0.3573  & 298 - 1100\\ 
         \ce{Ar}$^a$ & 20.79 & 0 & 0  & 298 - 1500\\ 
         \ce{CO}$^a$ & 26.87& 6.939  & -0.08237  & 298 - 1500\\ 
         \ce{C_{sus}}$^a$& -0.4493& 35.53 & -1.308 & 298 - 1500\\ 
         \ce{H_2O}$^a$&30.89 & 7.858  &0.2494 & 298 - 1300\\ 
         \ce{H_2}$^a$&   28.95& -0.5839&  0.1888 & 298 - 1500\\ \hline         
    \end{tabular}
    
    \footnotesize{ $^a$ based on data from  \cite{CRC2022}, $^b$ coefficients from \cite{HANEIN2020106043} }
    \label{tab:heatCap}
\end{table}

\end{document}